\documentclass[12pt]{amsart}
\usepackage{amssymb}
\newcommand{\R}{\mathbb{R}}
\newcommand{\N}{\mathbb{N}}
\newcommand{\Z}{\mathbb{Z}}
\newcommand{\C}{\mathbb{C}}

\newcommand{\F}{\mathcal{F}}
\newcommand{\f}{\mathcal{F}_{\alpha,q}}
\newcommand{\Ji}{J_{\alpha}^{(3)}}
\newcommand{\ji}{j_{\alpha}^{(3)}}

\def\endproof{\quad \hfill$\blacksquare$\vspace{0.15cm}\\}
\newcommand{\ds}{\displaystyle}

\newtheorem{defin}{Definition}
\newtheorem{propo}{ Proposition}
\newtheorem{lemme}{Lemma}
\newtheorem{coro}{Corollary}
\newtheorem{theorem}{Theorem}

\begin{document}

\title [Wavelet Transforms Associated With the Basic Bessel Operator ] {Wavelet Transforms Associated With the Basic Bessel Operator}
\author{Ahmed Fitouhi \quad \&\quad N\'eji Bettaibi  \quad \&\quad Wafa Binous  }
\address{A. Fitouhi. Facult\'e des Sciences de Tunis,
 1060 Tunis, Tunisia.}
\email{Ahmed.Fitouhi@fst.rnu.tn}
\address{N. Bettaibi. Institut
Pr\'eparatoire  aux \'Etudes d'Ing\'enieur de Nabeul, 8000 Nabeul,
Tunisia.}
 \email{Neji.Bettaibi@ipein.rnu.tn}
\address{W. Binous, Institut
Pr\'eparatoire  aux \'Etudes d'Ing\'enieur de Tunis, Tunis,
Tunisia.}

\begin{abstract}
\noindent This paper aims  to study the  $q$-wavelet and the
$q$-wavelet transforms, associated with the $q$-Bessel operator
for a fixed $q\in ]0, 1[$. As an application, an inversion
formulas of the $q$-Riemann-Liouville and $q$-Weyl transforms
using $q$-wavelets are given. For this purpose, we shall attempt
to extend the classical theory by giving their $q$-analogues.
 \end{abstract}
 \maketitle

\section {Introduction}\ \vspace{0.15cm}
Continuous wavelet transforms have been introduced by A. Grossmann
and J. Morlet \cite{GM} in the beginning of the 1980's and became
an active field of research, due to the fact that applications of
wavelet analysis to the diverse subjects of communication, seismic
data, signal and image processing... are being uncovered.\\

In \cite{Ft}, A. Fitouhi and K. Trim\`eche generalized the theory
as presented by T. H. Koornwinder \cite{koornwinder} and studied
the generalized wavelets and the generalized continuous wavelet
transforms associated with a class of singular differential
operators. This class contains, in particular, the so called
Bessel operator, which was studied extensively by K. Trim\`eche in
\cite{Tri}.\\

In this paper, we shall try to generalize our results in \cite{FN}
by studying   wavelets and continuous wavelet transforms
associated with the $q$-Bessel operator, studied in \cite{Fih}.
The basic tool in this work is some elements of $q$-harmonic
analysis related to the just mentioned operator. Next, using the
$q$-Riemann-Liouville and the $q$-Weyl operators, we will give
some relations between the continuous $q$-wavelet transform,
studied in \cite{FN}, and the continuous $q$-wavelet transform
associated with the $q$-Bessel operator, and  we deduce other
formulas which give the inverse operators of the
$q$-Riemann-Liouville and the $q$-Weyl transforms. These formulas
are better than those given in \cite{Fih} and \cite{fb} because
they are simple and we have a large choice of $q$-wavelets
associated with the $q$-Bessel operator, that  can be used in
these formulas.

 We are not in a situation to claim that all our results are new, but the methods
used are direct and constructive, and have a good resemblance with
the classical ones. Our approach in this paper is very similar to
the classical picture developed in \cite{Ft} and \cite{Tri}.\\

 This paper is organized as follows: in Section 2, we present some
 preliminaries results and notations that will be useful in the
 sequel. In Section 3, we establish some $q$-harmonic results
 associated with the $q$-Bessel operator.

 In Section 4, we  define the $q$-wavelets and the $q$-wavelet transforms associated
 with the $q$-Bessel operator, and discuss
 their properties. Special attention is paid  to the $q$-analogues of
 the Plancherel formula and the Parseval formula, and an inversion
 formula is proved. In Section 5, we give a characterization of the image set
 of the $q$-wavelet transform associated  with the $q$-Bessel operator.
 Section 6, is devoted to give some inversion formulas of the $q$-Riemann-Liouville and the $q$-Weyl
transforms. Finally, in Section 7, we give some relations between
 the continuous $q$-wavelet transform and the continuous $q$-wavelet transform associated
with the $q$-Bessel operator. We use these relations to derive the
inversion formulas of the $q$-Riemann-Liouville and the $q$-Weyl
transforms using wavelets.

  \section { Notation and preliminaries}
 \indent Throughout this paper, we will fix $q \in ]0,1[$ such that $\frac{Log(1-q)}{Log q}\in \Z $
 and $\alpha>-\frac{1}{2}$.
 We recall some usual notions and notations used in the
$q$-theory (see \cite {GR}).\\
For $a\in \C$, the $q$-shifted factorials are defined by
\begin {equation}
(a;q)_0=1;~~ (a;q)_n={\prod_{k=0}^{n-1}(1-aq^k)}, ~~
n=1,2,\ldots;~~{(a;q)_\infty} ={\prod_{k=0}^{\infty}(1-aq^k)}.
\end {equation}
 We also denote
\begin {equation}
  (a_1,a_2,\ldots ,a_p;q)_n=(a_1;q)_n(a_2;q)_n\ldots(a_p;q)_n , ~~~~   n={0},{1},{2},{3},\ldots
  \infty,
\end {equation}

\begin {equation}
 [x]_q={{1-q^x}\over{1-q}},~~ x\in \C ~~{\rm and} ~~[n]_q! ={{(q;q)_n}\over
 {(1-q)^n}},  ~~~~~~ n\in \N.
\end {equation}
The $q$-derivative $D_qf$  of a function $f$ is given by
\begin {equation}\label {21}
(D_qf)(x)={{f(x)-f(qx)}\over{(1-q)x}},~~  {\rm if}~~ x\not=0,
\end {equation}
$(D_qf)(0)=f'(0)$  provided $f'(0)$ exists. If  $f$ is
differentiable then $(D_qf)(x)$  tends to $f'(x)$ as $q$ tends to
1.\\
 The $q$-Jackson integrals from $0$ to $a$ and from $0$ to $\infty$ are
defined by (see \cite {Jac})
\begin {equation} \label {26}
\int_0^{a}{f(x)d_qx} =(1-q)a\sum_{n=0}^{\infty}{f(aq^n)q^n},
\end {equation}
\begin {equation} \label{int}
\int_0^{\infty}{f(x)d_qx}
=(1-q)\sum_{n=-\infty}^{\infty}{f(q^n)q^n},
\end {equation}
provided the sums converge absolutely.\\The $q$-Jackson integral
in a generic interval $[a,b]$ is given  by (see \cite {Jac})
\begin {equation}
\int_a^{b}{f(x)d_qx} =\int_0^{b}{f(x)d_qx}-\int_0^{a}{f(x)d_qx}.
\end {equation}
 \indent Jackson \cite {Jac} defined a $q$-analogue of the Gamma
function by
\begin {equation}
\Gamma_q (x) ={(q;q)_{\infty}\over{(q^x;q)_{\infty}}}(1-q)^{1-x} ,
\qquad x\not={0},{-1},{-2},\ldots.
\end {equation}
It is well known that it satisfies
\begin {equation}
\Gamma_q(x+1)=\frac{1-q^x}{1-q} \Gamma_q(x),\quad
\Gamma_q(1)=1~~~~{\rm and}~~~~
\lim_{q\rightarrow1^-}\Gamma_q(x)=\Gamma(x),~~ \Re(x)>0.
\end {equation}

We denote by
\begin {equation}
\R_q = \{\pm q^n: n\in \Z\}\cup\{0\},~~~~ \R_{q,+} = \{q^n: n\in
\Z\} ~~~~{\rm and }~~~~ \widetilde{\R}_{q,+}=\R_{q,+}\cup \{0\}.
\end {equation}
  $\bullet$ $\mathcal{E}_{*q}(\R_q)$ the space of the restrictions  on $\R_q$ of even infinitely
$q$-differentiable functions on $\R$, equipped with the induced
 topology  of  uniform convergence on all compact, for all functions and  its
   $q$-derivatives.\\
$\bullet$ $\mathcal{D}_{*q}(\R_q)$ the space of the restrictions
on $\R_q$  of even infinitely $q$-differentiable functions on $\R$
with compact supports,  equipped with the induced topology of
uniform convergence, for all
functions and its $q$-derivatives.\\
 $\bullet$ $\mathcal{C}_{*q,0}(\R_q)$ the space of the restrictions on
$\R_q$ of even smooth functions, continued in $0$ and vanishing at
$\infty$, equipped with  the induced topology of
uniform convergence. \\
$\bullet$ $\mathcal{S}_{*q}(\R_q)$ the space of the restrictions
on $\R_q$ of infinitely $q$-differentiable, even and fast
decreasing functions  and all its $q$-derivatives \textbf{i.e.}
$$\forall n,m\in\N,~~~~P_{n,m,q}(f)=\sup_{x\in\R; 0\leq
k\leq n}\mid (1+x^2)^m D_q^kf(x)\mid<+\infty.$$
$\mathcal{S}_{*q}(\R_q)$ is equipped with the induced topology
defined by the semi-norms $P_{n,m,q}$. \\
  $\bullet$  $L_q^p({\R}_{q,+},x^{2\alpha+1}d_qx)$, $p>0$, the set of all
functions defined on $\R_{q,+}$ such that
\begin {equation}
\|f\|_{p,\alpha,q}=\left\{\int_0^{\infty}\mid f(x)\mid^p
x^{2\alpha+1}d_qx\right\}^{\frac{1}{p}} < \infty.
\end {equation}
\section{Preliminaries on $q$-Harmonic Analysis Related to the $q$-Bessel Operator}
\subsection{Normalized $q$-Bessel function}
The $q$-Bessel operator is defined and studied in \cite{Fih} by
\begin {eqnarray*}
\Delta _{q,\alpha }f(z)&=&\left( \frac{1}{x^{2\alpha
+1}}D_{q}[x^{2\alpha +1}D_{q}f]\right) \left( q^{-1}z\right)\\
&=&q^{2\alpha + 1}\Delta _qf(z)+\frac{1-q^{2\alpha +
1}}{(1-q)q^{-1}z}D_qf(q^{-1}z),
\end{eqnarray*}
where
\begin {eqnarray*}
\Delta _qf(z)=D_q^2f(q^{-1}z).
\end{eqnarray*}
 We recall (see \cite{Fih}) that for $\lambda~\in \C$, the problem
\begin {equation}
         \left\{%
\begin{array}{ll}
   \Delta_{\alpha,q}u(x)  =  -\lambda^2u(x), &  \\
    u(0)=1, \hbox{$u'(0)=0$} \\
\end{array}%
\right.
\end{equation}
has as unique solution the  normalized q-Bessel function, given by
\begin {equation}
j_{\alpha }^{\left( 3\right) }(z;q^{2})=\left( 1-q^2\right)
^{\alpha }\Gamma _{q^{2}}\left( \alpha +1\right) \left((1-q)
q^{-1}z\right) ^{-\alpha }J_{\alpha }^{\left( 3\right) }\left(
\left( 1-q\right) q^{-1}z;q^{2}\right),
\end{equation}
 where
\begin{equation*}
\Ji\left( z;q^{2}\right) =\frac{z^{\alpha }\left(
q^{2\alpha +2};q^{2}\right) _{\infty }}{\left( q^{2};q^{2}\right) _{\infty }}%
\ _{1}\varphi _{1}\left( 0;q^{2\alpha +2};q^{2},q^{2}z^{2}\right)
\end{equation*}
is the  Jackson's third q-Bessel function. This  function is
called in some literature the Hahn-Exton q-Bessel function ( see
\cite{KS}).

The following lemma shows some estimations for the normalized
$q$-Bessel function.
\begin{lemme}\label{estim} For $x \in \R_{q,+}$, we have\\
1)   $\ds |j_\alpha^{(3)}(x;q^2)|\leq \frac{1}{(q;q^2)_\infty^2}$;\\
2)   $\ds |j_\alpha^{(3)}(x;q^2)|\leq
\frac{(-q^2;q^2)_\infty(-q^{2(\alpha+1)};q^2)_\infty}{(q^{2(\alpha+1)};q^2)_\infty}\left\{%
\begin{array}{ll}
    1, & \hbox{if ~~$ x\leq \frac{q}{1-q}$ ,} \\
    q^{\left( \frac{Log(\frac{1-q}{q}x)}{Log q}\right)^2}, & \hbox{if ~~$ x\geq \frac{q}{1-q}$;} \\
\end{array}%
\right.$\\
3) For all $\nu\in\R$, we have $\ds j_\alpha^{(3)}(x;q^2)=
o(x^{-\nu})$ as $x\rightarrow +\infty$.\\
In particular, we have $\ds \lim_{x\rightarrow
+\infty}j_\alpha^{(3)}(x;q^2)=0$.
\end{lemme}
\proof.\\
1) is proved in \cite{Fih}.\\
2) From the properties of the basic function $_1\varphi_1$ ( see
\cite{Fih} or  \cite{KS}), we have:\\
 $\centerdot$ For $x =q^n \in \R_{q,+}$ $n\in\N$,
\begin{eqnarray*}
|x^{-\alpha}J_\alpha^{(3)}(x;q^2)|&=& \frac{1}{(q^2;q^2)_\infty}|
(q^{2\alpha+2};q^2)_\infty~~
_1\varphi_1(0; q^{2\alpha+2};q^2,q^{2n+2})|\\
&\leq&
\frac{1}{(q^2;q^2)_\infty}(-q^{2(n+1)};q^2)_\infty(-q^{2\alpha+2};q^2)_\infty\\
&\leq&
\frac{1}{(q^2;q^2)_\infty}(-q^2;q^2)_\infty(-q^{2\alpha+2};q^2)_\infty.
\end{eqnarray*}
$\centerdot$ For $x =q^{-n} \in \R_{q,+}$ $n\in\N$,
\begin{eqnarray*}
|x^{-\alpha}J_\alpha^{(3)}(x;q^2)|&=& \frac{1}{(q^2;q^2)_\infty}|
(q^{2(1-n)};q^2)_\infty~~
_1\varphi_1(0;q^{2(1-n)} ;q^2,q^{2\alpha+2})|\\
&\leq& \frac{1}{(q^2;q^2)_\infty}
q^{n(n+2\alpha+1)}(-q^2;q^2)_\infty(-q^{2\alpha+2};q^2)_\infty\\
&\leq& \frac{1}{(q^2;q^2)_\infty}
q^{n^2}(-q^2;q^2)_\infty(-q^{2\alpha+2};q^2)_\infty,
\end{eqnarray*}
since $\alpha>-1/2$.

So, $$|x^{-\alpha}J_\alpha^{(3)}(x;q^2)|\leq
\frac{(-q^2;q^2)_\infty(-q^{2(\alpha+1)};q^2)_\infty}{(q^2;q^2)_\infty}\left\{
\begin{array}{ll}
    1, & \hbox{if ~~$n\geq 0$,} \\
    q^{n^2}, & \hbox{if ~~$n\leq 0$,} \\
\end{array}%
\right.
$$
which is equivalent to
$$|x^{-\alpha}J_\alpha^{(3)}(x;q^2)|\leq
\frac{(-q^2;q^2)_\infty(-q^{2(\alpha+1)};q^2)_\infty}{(q^2;q^2)_\infty}\left\{
\begin{array}{ll}
    1, & \hbox{if ~~$x\leq 1$,} \\
    q^{\left( \frac{Log(x)}{Log q}\right)^2}, & \hbox{if ~~$x\geq 1$.} \\
\end{array}%
\right.
$$
The relation 2)  follows from this inequality and the relation:
$$j_\alpha^{(3)}(z;q^2)
=(1-q^2)^\alpha\Gamma_{q^2}(\alpha+1)(\frac{1-q}{q}z)^{-\alpha}J_\alpha^{(3)}(\frac{1-q}{q}z;q^2).$$
Relation 3) is a direct consequence of 2).
\endproof

\begin{lemme}\label{ort} For $x~~,y\in\R_{q,+}$, we have
\begin{equation}
(xy)^{\alpha +1}\int_0^\infty \ji (xt,q^2)\ji
(yt,q^2)t^{2\alpha+1}d_qt=
\frac{(1+q)^{2\alpha}\Gamma_{q^2}^2(\alpha+1)q^{2(\alpha+1)}}{1-q}\delta_{x,y}.
\end{equation}
\end{lemme}
\proof The result follows from the definition of $\ji$ and the
orthogonality relation of $\Ji$ proved in \cite{KS}.
\endproof
\subsection{$q$-Bessel Fourier transform}

The generalized $q$-Bessel translation operator
$T_{q,x}^\alpha,~~x\in \R_{q,+}$ was defined
 in \cite {Fih} on $\mathcal{D}_{*q}(\R_q)$ by
\begin {equation}\label{trq}
T_{q,x}^\alpha(f)(y)=\sum_{n=0}^{\infty} \frac{q^{n^2}}{(q^2,
q^{2\alpha+2};q^2)_n}
(\frac{x}{y})^{2n}\sum_{k=-n0}^{n}(-1)^{n-k}U_k(n) f(q^ky), y\in
\R_{q,+}
\end {equation}
and $T_{q,0}^\alpha(f)=f$, where
$$U_k(n)=q^{k(k-1)+2n(k+\alpha)}\sum_{p=0}^k\left[ \begin{array}{c}
  n \\
  p
\end{array}\right]_{q^2}\left[ \begin{array}{c}
  n \\
  n+k-p
\end{array}\right]_{q^2}q^{-2p(n+k+\alpha-p)}$$
is the  $q$-Bessel $q$-Binomial coefficient associated with the
$q$-Bessel operator (see \cite{Fih}). It verifies, in particular

\begin {equation}
\int_0^\infty
T_{q,x}^\alpha(f)(y)g(y)y^{2\alpha+1}d_qy=\int_0^{\infty}
f(y)T_{q,x}^\alpha(g)(y)y^{2\alpha+1}d_qy,~~
x\in\widetilde{\R}_{q,+},
\end {equation}
and
\begin {equation}
T_{q,x}^\alpha\ji(ty;q^2)=\ji(tx;q^2)\ji(ty;q^2),~~ x, y,
t\in\widetilde{\R}_{q,+}.
\end {equation}
\indent The $q$-Bessel Fourier transform and the $q$-convolution
product are defined \\ (see \cite {Fih}) for $f,g \in
\mathcal{D}_{*q}(\R_q)$, by
\begin {equation}
\f(f)(\lambda)=c_{\alpha,q}\int_0^\infty f(x)\ji(\lambda
x;q^2)x^{2\alpha+1}d_qx,
\end {equation}
\begin {equation}
f*_Bg(x)=c_{\alpha,q}\int_0^\infty T_{q,x}^\alpha
f(y)g(y)y^{2\alpha+1}d_qy,
\end {equation}
where
\begin {equation}\label{cq}
c_{\alpha,q}=\frac{(1+q^{-1})^{-\alpha}}{\Gamma_{q^2}(\alpha+1)}.
\end {equation}
Using the proprieties of the $q$-generalized Bessel translation,
one can prove easily the following result \cite{Fih}.

\begin{theorem} For $f, g\in\mathcal{D}_{*q}(\R_q)$, we have
\begin {equation}
\f(f*_Bg)=\f(f)\f(g),
\end {equation}
\begin {equation}\label{tran}
\f(T_{q,x}^\alpha f)(\lambda)=\ji(\lambda x;q^2)\f(f)(\lambda),
~~~~x\in \widetilde{\R}_{q,+},~~\lambda \in \R_{q,+}
\end {equation}
and
\begin {equation}
\f(\Delta_{\alpha,q}f)(\lambda)=-\frac{\lambda^2}{q^{2\alpha+1}}\f(f)(\lambda),~~~~\lambda~\in\C.
\end {equation}
\end {theorem}

\begin{theorem} For $f\in L_q^1(\R_{q,+},x^{2\alpha+1}d_qx)$, we
have\\
 \begin {equation} \label{l1}
  \F_{\alpha,q}(f)\in\mathcal{C}_{*q,0}(\R_q)
\end {equation}
and
\begin {equation}
  \|\F_{\alpha,q}(f)\|_{\mathcal{C}_{*q,0}(\R_q)}\leq
  \frac{c_{\alpha,q}}{(q;q^2)_\infty^2}\|f\|_{1,\alpha, q}.
\end {equation}
\end {theorem}
\proof Let $f\in L_q^1(\R_{q,+},x^{2\alpha+1}d_qx)$. From the
relation 1) of Lemma \ref{estim}, we have
$$\forall\lambda, ~~x\in \R_{q,+},~~~~  |f(x)\ji (\lambda x;q^2)
x^{2\alpha +1}|\leq \frac{1}{(q;q^2)_\infty^2}|f(x) x^{2\alpha
+1}|.$$
 Then, the definition of $\ji$, the relation 3)
of Lemma \ref{estim} and the Lebesgue theorem imply that $\f
(f)\in\mathcal{C}_{*q,0}(\R_q)$. On the other hand, we have
  for all $\lambda\in \R_{q,+}$,
$$|\f (f)(\lambda)|\leq \frac{c_{\alpha,q}}{(q;q^2)_\infty^2}\|f\|_{1,\alpha,
q},$$  which achieves the proof.
\endproof

\begin{theorem}\label{ld}.\\
 $\f$ is an isomorphism of $L_q^2(\R_{q,+}, x^{2\alpha+1}d_qx)$
(resp. $\mathcal{S}_{*q}(\R_q)$), $\f^{-1}=q^{-4\alpha-2}\f$ and
for $f\in L_q^2(\R_{q,+}, x^{2\alpha+1}d_qx)$, we have
\begin {equation}
  \|\f(f)\|_{2,\alpha,q}=q^{2\alpha+1}\|f\|_{2,\alpha, q}.
\end {equation}
\end {theorem}
\proof.\\
The  parity of $\ji$ and the relation
$$\f (\Delta_{\alpha, q}f)=-\frac{\lambda^2}{q^{2\alpha+1}}\f(f)$$
show that if $f$ is in $\mathcal{S}_{*q}(\R_q)$, then $\f (f)$
belongs to $\mathcal{S}_{*q}(\R_q)$.\\
 Lemma \ref{ort} achieves the proof.\endproof
{\bf Remak 1.}\\
Using the previous theorem and the relation (\ref{tran}), one can
see that, for $f\in L_q^2(\R_{q,+}, x^{2\alpha+1}d_qx)$ (resp.
$\mathcal{S}_{*q}(\R_q)$), we have for all $x\in
\widetilde{\R}_{q,+}$, $T_{q,x}^\alpha f$ belongs to
$L_q^2(\R_{q,+}, x^{2\alpha+1}d_qx)$ (resp.
$\mathcal{S}_{*q}(\R_q)$) and
\begin {equation}\label{tcos}
  \|T_{q,x}^\alpha f\|_{2,\alpha,q}\leq\frac{1}{(q;q^2)_\infty^2} \|f\|_{2,\alpha,q}.
\end {equation}

\begin{propo}
Let $f$ and $g$ be in $L_q^2(\R_{q,+}, x^{2\alpha +1}d_qx)$, then\\
 1)  $f*_Bg \in L_q^2(\R_{q,+}, x^{2\alpha+1}d_qx)$   iff   $\f(f)\f(g) \in
 L_q^2(\R_{q,+}, x^{2\alpha+1}d_qx)$.\\
 2)
\begin {equation}\label{prop2}
q^{4\alpha+2}\int_0^\infty \mid f*_Bg(x)
\mid^2x^{2\alpha+1}d_qx=\int_0^\infty \mid \f(f)(x)\mid^2 \mid
\f(g)(x) \mid^2 x^{2\alpha+1} d_qx,
\end {equation}
where both sides are finite or infinite.
\end{propo}
\proof The proof is a direct consequence of Theorem \ref{ld} and
the fact that \\ $\f (f*_Bg)=\f (f)\f (f)$. \endproof

\section {$q$-Wavelet transforms associated with the $q$-Bessel opertor}
\begin {defin}
A $q$-wavelet associated with the $q$-Bessel operator is an even
function $g\in  L_q^2(\R_{q,+}, x^{2\alpha+1}d_qx)$ satisfying the
following admissibility condition:
\begin{equation}\label{wc}
0<C_g=\int_0^{\infty}\mid \f(g)(a)\mid^2 \frac{d_qa}{a}<\infty.
\end{equation}
 \end {defin}
 {\bf Remarks}\\
  1) For all $\lambda \in \R_{q,+}$, we have
 $$ C_g=\int_0^{\infty}\mid \f(g)(a\lambda)\mid^2
\frac{d_qa}{a}.$$
 2) Let $f$ be a nonzero function in $\mathcal{S}_{*q}(\R_q)$
 (resp. $\mathcal{D}_{*q}(\R_q)$). Then $g=\Delta_{\alpha,q}f$ is a
 $q$-wavelet associated with the $q$-Bessel operator, in $\mathcal{S}_{*q}(\R_q)$
 (resp. $\mathcal{D}_{*q}(\R_q)$) and we have
 $$C_g= \frac{1}{q^{4\alpha+2}}\int_0^\infty a^3\mid\f(f)(a) \mid^2d_qa .$$

{\bf Example} \\
Consider the functions  $ \ds G(x;q^2)= A_\alpha
e_{q^2}^{-\frac{q^{-(2\alpha+1)}}{(1+q)^2}x^2}$ and $g=
\Delta_{\alpha,q} G(.;q^2)$ , \ where\\ $\ds
A_\alpha=c_{\alpha,q}\int_0^\infty
x^{2\alpha+1}e_{q^2}^{-x^2}d_qx$
and $e_{q^2}$ is the $q$-analogue of the exponential function.\\
We have   $x\mapsto G(x;q^2)$ is in $\mathcal{S}_{*q}(\R_q)$ and
(see \cite {Fih}, Proposition 8)
$$ \f(G(.;q^2))(x)=q^{4\alpha+2}e_{q^2}^{-x^2}, ~~~~~~x\in \R_{q,+}.$$
Then,  $g$ is in $\mathcal{S}_{*q}(\R_q)$ and
$$\f(g)(x)=-\frac{x^2}{q^{2\alpha+1}}\f(G(.;q^2))(x)=-q^{2\alpha+1}x^2e_{q^2}^{-x^2},~~~~~~x\in \R_{q,+}.$$
It is then easy to see that
$$0< \mid \f(g)\mid^2(a)\leq q^{4\alpha+2} a^4e_{q^2}^{-a^2}, ~~~~~~\forall a
\in \R_{q,+}.$$

Thus
\begin {eqnarray*}
0< \int_0^{\infty}\mid \f(g)\mid^2(a) \frac{d_qa}{a}&\leq&
q^{4\alpha+2}
\int_0^{\infty}a^3e_{q^2}^{-a^2}d_qa \\
&=&
\frac{q^{4\alpha+2}}{(1+q)}\frac{(-q^4,-q^{-2};q^2)_{\infty}}{(-q^2,-1;q^2)_{\infty}}\\
&=&\frac{q^{4\alpha}}{(1+q)}.
\end{eqnarray*}
So $g$ is a $q$-wavelet associated with the $q$-Bessel operator.

\begin{propo}\label{ne}
Let $g\neq 0$ be a function in $\L_q^2(\R_{q,+}, x^{2\alpha
+1}d_qx)$ satisfying:\\
    1) $\f(g)$ is continuous at 0.\\
    2) $\exists\beta>0$ such that
    $\f(g)(x)-\f(g)(0)=O(x^{\beta})$, as $x\rightarrow 0$.\\
Then,  (\ref{wc}) is equivalent to
\begin {equation}\label{f0}
 \f(g)(0)=0.
\end {equation}
\end{propo}
\proof $\centerdot$ We suppose that (\ref{wc}) is satisfied.\\
If $\f(g)(0)\neq0$, then from the condition  1) there exist
$p_0\in\N$ and $M>0$, such that
$$ \forall n\geq p_0, ~~\mid \f(g)(q^n)\mid \geq M.$$
Then, the integral in (\ref{wc}) would be equal to $\infty$.\\
$\centerdot$ Conversely, we suppose that $\f(g)(0)=0$.\\
As $g\neq0$, we deduce from Theorem \ref{ld}, that the first
inequality in (\ref{wc}) is satisfied.\\
 On the other hand, from the condition 2), there exist $n_0\in\N$ and
 $\epsilon>0$, such that for all $n\geq n_0$,
 $$ \mid \f(g)(q^n)\mid \leq \epsilon q^{n\beta}.$$
Then using the definition of the $q$-integral and Theorem
\ref{ld}, we obtain
\begin {eqnarray*}
\int_0^\infty \mid
\f(g)(a)\mid^2\frac{d_qa}{a}&=&(1-q)\sum_{n=-\infty}^{\infty}\mid
\f(g)(q^n)\mid^2\\
&=& (1-q)\sum_{n=-\infty}^{n_0}\mid
\f(g)(q^n)\mid^2+(1-q)\sum_{n=n_0+1}^{\infty}\mid
\f(g)(q^n)\mid^2\\
&\leq&
\frac{(1-q)}{q^{(2\alpha+2)n_0}}\sum_{n=-\infty}^{\infty}q^{(2\alpha+2)n}\mid
\f(g)(q^n)\mid^2+(1-q)\epsilon^2\sum_{n=0}^{\infty}q^{2n\beta}\\
&\leq&
\frac{\|\f(g)\|_{2,\alpha,q}^2}{q^{(2\alpha+2)n_0}}+\frac{1-q}{1-q^{2\beta}}\epsilon^2\\
&=&q^{(4\alpha+2)}\frac{\|g\|_{2,\alpha,q}^2}{q^{(2\alpha+2)n_0}}+\frac{1-q}{1-q^{2\beta}}\epsilon^2.
\end {eqnarray*}
This proves the second inequality of (\ref{wc}).\endproof

{\bf Remark 2.}\\ Owing to (\ref{l1}), the continuity assumption
in the previous proposition will certainly hold if $g$ is moreover
in
 $L_q^1(\R_{q,+},
x^{2\alpha+1}d_qx)$. Then (\ref{f0}) can be equivalently written
as
$$ \int_0^{\infty}g(x)x^{2\alpha+1}d_qx=0.$$
\begin{theorem}
Let $a\in \R_{q,+}$ and $ g\in  L_q^2(\R_{q,+},
x^{2\alpha+1}d_qx)$. Then, the function $g_a$ defined for $x\in
\R_{q,+}$, by
\begin{equation}\label{ga}
g_a(x)=\frac{1}{a^{2\alpha+2}}g(\frac{x}{a}),
\end{equation}
satisfies:\\
i) the function $g_a$ belongs to $ L_q^2(\R_{q,+},
x^{2\alpha+1}d_qx)$  and we have
\begin {equation}\label{g-a}
\| g_a \|_{2,\alpha,q}=\frac{1}{a^{\alpha+1}}\| g \|_{2,\alpha,q};
\end{equation}
ii)  for all $\lambda\in \R_{q,+}$, we have
\begin {equation} \label{pr1}
\f (g_a)(\lambda)=\f (g)(a\lambda).
\end{equation}
\end{theorem}
\proof The change of variable $\ds u=\frac{x}{a}$ leads to:
\begin{eqnarray*}
\int_0^\infty |g_a(x)|^2x^{2\alpha+1}d_qx&=&
\frac{1}{a^{4\alpha+4}}\int_0^\infty
|g(\frac{x}{a})|^2x^{2\alpha+1}d_qx\\
&=&\frac{1}{a^{2\alpha+2}}\int_0^\infty |g(u)|^2u^{2\alpha+1}d_qu
\end{eqnarray*}
and for $\lambda\in \R_{q,+}$,
\begin{eqnarray*}
\f(g_a)(\lambda)&=&\frac{c_{\alpha,q}}{a^{2\alpha+2}}\int_0^\infty
g(\frac{x}{a})\ji(\lambda x;q^2)x^{2\alpha+1}d_qx\\
&=&c_{\alpha,q}\int_0^\infty g(u)\ji(a\lambda
u;q^2)u^{2\alpha+1}d_qu = \f(g)(a\lambda).
\end{eqnarray*}
\endproof
\begin{propo}
Let $g$ be in $\mathcal{S}_{*q}(\R_q)$ (resp
 $\mathcal{D}_{*q}(\R_q)$). Then for all $a\in \R_{q,+}$ the
 function $g_a$ given by the relation (\ref{ga}) belongs to $\mathcal{S}_{*q}(\R_q)$ (resp
 $\mathcal{D}_{*q}(\R_q)$).
\end{propo}
\begin{theorem}
Let $g$ be a $q$-wavelet associated with the $q$-Bessel operator
in  $L_q^2(\R_{q,+}, x^{2\alpha+1}d_qx)$
(resp.$\mathcal{S}_{*q}(\R_q)$). Then for all $a\in \R_{q,+}$ and
$b\in \widetilde{\R}_{q,+}$, the function
\begin{equation}\label{wavelet}
g_{a,b}(x)=\sqrt{a}T_{q,b}^\alpha(g_a),
\end{equation}
is a $q$-wavelet associated with the $q$-Bessel operator in
$L_q^2(\R_{q,+}, x^{2\alpha+1}d_qx)$ (resp.
$\mathcal{S}_{*q}(\R_q)$) and we have
\begin {equation}\label{wab}
C_{g_{a,b}}= a\int_0^\infty
(\ji(\frac{xb}{a};q^2))^2\mid\f(g)(x)\mid^2\frac{d_qx}{x}.
\end{equation}
Where $T_{q,b}^\alpha, ~~b\in \widetilde{\R}_{q,+}$ are the
$q$-generalized translations defined by the relation (\ref{trq}).
\end{theorem}
\proof As $g_a$ is in $L_q^2(\R_{q,+}, x^{2\alpha+1})$ (resp.
$\mathcal{S}_{*q}(\R_q)$), Remark 1 shows that the relation
(\ref{wavelet}) defines an element of $L_q^2(\R_{q,+},
x^{2\alpha+1})$ (resp. $\mathcal{S}_{*q}(\R_q)$). On the other
hand, from the relations (\ref{pr1}) and  (\ref{tran}), we have
for all $\lambda\in \R_{q,+}$,
$$\f(g_{a,b})(\lambda)=\sqrt{a}\ji(b\lambda;q^2)\f(g)(a\lambda).$$
This relation implies (\ref{wab}).\\
Now, we shall prove that the function $g_{a,b}$ satisfies the
admissibility relation (\ref{wc}). \\
As $g\neq 0$, we deduce from (\ref{wab}) and Theorem \ref{ld} that
$C_{g_{a,b}}\neq 0$. On the other hand, from the relation
(\ref{wc}) and the relation 1) of Lemma \ref{estim}, we deduce
that
$$ C_{g_{a,b}} \leq \frac{a}{(q;q^2)_\infty^4}C_g,$$
which gives the result.\endproof

\begin{propo}\label{cont}
Let $g$ be a $q$-wavelet associated with the $q$-Bessel operator
in $L_q^2(\R_{q,+}, x^{2\alpha +1}d_qx)$. Then the mapping
$$F:(a,b)\mapsto g_{a,b}$$ is continuous from
$\R_{q,+}\times\widetilde{\R}_{q,+}$ into $L_q^2(\R_{q,+},
x^{2\alpha +1}d_qx)$.
\end{propo}
\begin{proof}It is clear that $F$ is a mapping from
$\R_{q,+}\times\widetilde{\R}_{q,+}$ into $L_q^2(\R_{q,+},
x^{2\alpha +1}d_qx)$ and it is continuous at all $\ds (a,b)\in
\R_{q,+}\times\R_{q,+}$.\\ Now,  fix $\ds a\in \R_{q,+}$. For
$b\in \widetilde{\R}_{q,+}$, we have
\begin{eqnarray*}
\parallel F(a,b)-F(a,0)\parallel_{2,\alpha, q}^2&=&\parallel
T_{q,b}^\alpha(g_a)-g_a\parallel_{2,\alpha, q}^2\\
&=&q^{-4\alpha-2}\parallel
\f\left(T_{q,b}^\alpha(g_a)-g_a\right)\parallel_{2,\alpha, q}^2\\
&=&q^{-4\alpha-2}\int_0^\infty \mid 1-\ji(xb;q^2)\mid^2 \mid
\f(g_a)\mid^2(x)x^{2\alpha+1}d_qx.
\end{eqnarray*}
However, for all $x\in \R_{q,+}$ and $b\in \widetilde{\R}_{q,+}$,
we have
$$\mid 1-\ji(xb;q^2)\mid^2 \mid
\f(g_a)\mid^2(x)\leq (1+\frac{1}{(q;q^2)_\infty ^2})^2\mid
\f(g_a)\mid^2(x)$$ and $\ds \f(g_a) \in
L_q^2(\R_{q,+},x^{2\alpha+1}d_qx)$. So, the Lebesgue theorem leads
to
$$\lim_{\small{\begin{array}{c}
          b\rightarrow 0 \\
          b\in \widetilde{\R}_{q,+}\\
        \end{array}}} \parallel F(a,b)-F(a,0)\parallel_{2,\alpha,q} =0.$$
Then for all open neighborhood $V$ of $F(a,0)$ in
$L_q^2(\R_{q,+},x^{2\alpha+1}d_qx)$, there exists an open
neighborhood $U$ of $0$ in $\widetilde{\R}_{q,+}$ such that $$
\forall b\in U, ~~F(a,b)\in V.$$ Thus $\{a\}\times U$ is an open
neighborhood of $(a,0)$ in $\R_{q,+}\times\widetilde{\R}_{q,+}$
and $F(\{a\}\times U)\subset V$. Which proves the continuity of
$F$ at $(a,0)$.
\end{proof}

\begin {defin}
Let $g$ be a $q$-wavelet associated with the $q$-Bessel operator
in $\mathcal{D}_{*q}(\R_q)$. We define the continuous $q$-wavelet
transform associated with the $q$-Bessel  operator  by
\begin {equation}\label{contin}
\Psi^{\alpha}_{q,g}(f)(a,b)=c_{\alpha,q} \int_0^\infty
f(x)\overline{g_{a,b}}(x)x^{2\alpha+1}d_qx, ~~~~~~a\in\R_{q,+},
~~b\in\widetilde{\R}_{q,+} ~~~~{\rm and} ~~~~~f \in
\mathcal{D}_{*q}(\R_q).
\end {equation}
\end {defin}
{\bf Remark 3.} The relation (\ref{contin}) can also be written in
the form
\begin {eqnarray*}
\Psi^{\alpha}_{q,g}(f)(a,b)&=& \sqrt{a}f*_B\overline{g_a}(b)\\
&=&\sqrt{a}q^{-4\alpha-2}\f(\f (f*_B\overline{g_a}))(b)\\
&=&\sqrt{a}q^{-4\alpha-2}\f\left[\f (f).\f (\overline{g_a})\right](b)\\
&=& \sqrt{a}~~q^{-4\alpha-2} c_{\alpha,q}\int_0^{\infty}\f
(f)(x).\f (\overline{g})(ax)\ji(bx;q^2)x^{2\alpha+1}d_qx,
\end {eqnarray*}
 where $c_{\alpha,q}$ is given by (\ref{cq}).

We give some properties of $\Psi^{\alpha}_{q,g}$ in the  following
proposition.
\begin{propo}
Let $g$ be a $q$-wavelet associated with the $q$-Bessel operator
in $L_q^2(\R_{q,+}, x^{2\alpha +1}d_qx)$ and $f\in
L_q^2(\R_{q,+}, x^{2\alpha +1}d_qx)$, then\\
i) For all $a\in \R_{q,+}$ and $b\in \widetilde{\R}_{q,+}$, we
have
\begin {equation}
\mid \Psi^{\alpha}_{q,g}(f)(a,b)\mid\leq
\frac{c_{\alpha,q}}{(q;q^2)_\infty^2a^{\alpha+1/2}}\|f\|_{2,\alpha,q}\|g\|_{2,\alpha,q}.
\end {equation}
ii) For all $a\in \R_{q,+}$, the function $b\mapsto
\Psi^{\alpha}_{q,g}(f)(a,b)$ is continuous on
$\widetilde{\R}_{q,+}$ and we have
\begin {equation}
\lim_{b\rightarrow\infty} \Psi^{\alpha}_{q,g}(f)(a,b)=0.
\end {equation}
 iii) If $g$ is in $\mathcal{S}_{*q}(\R_{q})$, then for all $f$ in
$\mathcal{S}_{*q}(\R_{q})$, the function $b\mapsto
\Psi^{\alpha}_{q,g}(f)(a,b)$ is in $\mathcal{S}_{*q}(\R_{q})$.
\end{propo}

\proof .\\
i) For $a\in \R_{q,+}$ and $b\in \widetilde{\R}_{q,+}$, we have
\begin {eqnarray*}
\mid \Psi^{\alpha}_{q,g}(f)(a,b)\mid &=& c_{\alpha,q} \mid
\int_0^\infty
f(x)\overline{g_{a,b}}(x)x^{2\alpha +1}d_qx \mid\\
&\leq& c_{\alpha,q} \sqrt{a}\int_0^\infty \mid f(x)\mid
\mid T_{q,b}g_a(x)\mid x^{2\alpha +1}d_qx\\
&\leq& \frac{c_{\alpha,q}}{(q;q^2)_\infty^2a^{\alpha+1/2}}
\|f\|_{2,\alpha,q}\|g\|_{2,\alpha,q},
\end {eqnarray*}
by using the relations (\ref{tcos}) and (\ref{g-a}).\\
ii) As in Proposition \ref{cont}, it suffices to prove the
continuity at $0$. For $b\in\widetilde{\R}_{q,+}$, we have
\begin {eqnarray*}
\Psi^{\alpha}_{q,g}(f)(a,b)&=&\sqrt{a}q^{-4\alpha-2}\f\left[\f
(f).\f (\overline{g_a})\right](b)\\
&=&\sqrt{a}q^{-4\alpha-2}c_{\alpha,q}\int_0^\infty \f (f)(x).\f
(\overline{g_a})(x)\ji (bx;q^2)x^{2\alpha+1}d_qx
\end {eqnarray*}
and
$$ \forall x\in \R_{q,+},~~~~\mid \ji(bx;q^2)\mid\leq\frac{1}{(q;q^2)_\infty
^2}.$$ Since $f,g\in L_q^2(\R_{q,+}, x^{2\alpha +1}d_qx)$, then by
Theorem \ref{ld},
 $\f(f)$ and $\f(\overline{g_a})$ are in $L_q^2(\R_{q,+}, x^{2\alpha +1}d_qx)$.\\
So, the product $\f(f).\f(\overline{g_a})$ is in $L_q^1(\R_{q,+},
x^{2\alpha +1}d_qx)$. Thus, by application of the Lebesgue
theorem, we obtain
\begin{eqnarray*}\lim_{\begin{array}{c}
        b\rightarrow 0 \\
          b\in\widetilde{\R}_{q,+} \\
        \end{array}} \Psi^{\alpha}_{q,g}(f)(a,b)&=&\lim_{\begin{array}{c}
        b\rightarrow 0 \\
          b\in\widetilde{\R}_{q,+} \\
        \end{array}} \sqrt{a}q^{-4\alpha-2}c_{\alpha,q}\int_0^\infty \f (f)(x).\f
(\overline{g_a})(x)\ji
(bx;q^2)x^{2\alpha+1}d_qx\\&=&\Psi^{\alpha}_{q,g}(f)(a,0).
 \end{eqnarray*}
 Which proves the continuity of $\Psi^{\alpha}_{q,g}(f)(a,.)$ at $0$.\\
 Finally   (\ref{l1}) implies that $$
\Psi^{\alpha}_{q,g}(a,b)= \sqrt{a}q^{-4\alpha-2}
\f[\f(f).\f(\overline{g_a})](b)$$ tends
to $0$ as $b$ tends to $\infty$.\\
iii)  is an immediate consequence of the relation
$$\Psi^{\alpha}_{q,g}(f)(a,b)= \sqrt{a}f*_B\overline{g_a}(b)$$
and the properties of the $q$-Bessel convolution product.
\endproof

\begin {theorem}\label{perplan}
Let $g\in L_q^2(\R_{q,+}, x^{2\alpha +1}d_qx)$ a $q$-wavelet associated with the $q$-Bessel operator.\\
i) \underline{Plancheral formula for $\Psi^{\alpha}_{q,g}$}\\
For $f\in L_q^2(\R_{q,+}, x^{2\alpha +1}d_qx)$, we have
\begin {equation}\label{plancherel}
\frac{1}{C_g}\int_0^\infty \int_0^\infty \mid
\Psi^{\alpha}_{q,g}(f)
(a,b)\mid^2b^{2\alpha+1}\frac{d_qbd_qa}{a^2}=\|f\|_{2,\alpha,
q}^2.
\end {equation}
ii)\underline{Parseval formula for $\Psi^{\alpha}_{q,g}$}\\
For $f_1, f_2\in L_q^2(\R_{q,+}, x^{2\alpha +1}d_qx)$, we have
\begin {equation}\label{Parseval}
\int_0^\infty f_1(x)\overline{f}_2(x)x^{2\alpha+1}d_qx=
\frac{1}{C_g}\int_0^\infty \int_0^\infty \Psi^{\alpha}_{q,g}(f_1)
(a,b)\overline{\Psi^{\alpha}_{q,g}}(f_2)
(a,b)b^{2\alpha+1}\frac{d_qad_qb}{a^2}.
\end {equation}
\end {theorem}
\proof By using Fubini's theorem, Theorem \ref{ld}, and the
relations (\ref {pr1}) and (\ref{prop2}), we have {\small \begin
{eqnarray*} q^{4\alpha+2}\int_0^\infty \int_0^\infty \mid
\Psi^{\alpha}_{q,g}(f)
(a,b)\mid^2&b^{2\alpha+1}&\frac{d_qad_qb}{a^2}=q^{4\alpha+2}\int_0^\infty
\left(\int_0^\infty
\mid f*_B\overline{g_a} \mid^2(b)b^{2\alpha+1}d_qb\right)\frac{d_qa}{a}\\
&=& \int_0^\infty \left(\int_0^\infty \mid
\f(f)(x)\mid^2\mid\f(\overline{g_a})
\mid^2(x)x^{2\alpha+1}d_qx\right)\frac{d_qa}{a}\\
&=& \int_0^\infty\mid \f(f)(x)\mid^2 \left(\int_0^\infty
\mid\f(g)(ax)
\mid^2\frac{d_qa}{a}\right)x^{2\alpha+1}d_qx\\
&=& C_g \int_0^\infty\mid \f(f)(x)\mid^2x^{2\alpha+1}d_qx= C_g
q^{4\alpha+2}\|f\|_{2,\alpha,q}^2.
\end {eqnarray*}}
The relation (\ref{plancherel}) is then proved.\\
ii) The result is easily deduced from (\ref{plancherel}).

{\bf Remark 4.}\\
If $g\in \L_q^2(\R_{q,+}, x^{2\alpha +1}d_qx)$ is a $q$-wavelet
associated with the $q$-Bessel operator, then for all $f\in
\L_q^2(\R_{q,+}, x^{2\alpha +1}d_qx)$, we have  $\ds
\Psi^{\alpha}_{q,g}(f)\in
\L_q^2(\R_{q,+}\times\widetilde{\R}_{q,+};
b^{2\alpha+1}\frac{d_qad_qb}{a^2})$ and
$$ \|\Psi^{\alpha}_{q,g}(f)\|_{\L_q^2(\R_{q,+}\times\widetilde{\R}_{q,+};
b^{2\alpha+1}\frac{d_qad_qb}{a^2})}^2=C_g\|f\|_{2,\alpha,q}^2.$$
\begin {theorem}\label{omri} Let $g$ be a $q$-wavelet associated with the $q$-Bessel operator in $\L_q^2(\R_{q,+}, x^{2\alpha +1}d_qx)$,
then for all $f\in \L_q^2(\R_{q,+}, x^{2\alpha +1}d_qx)$, we have
\begin {equation}
f(x)=\frac{c_{\alpha,q}}{C_g}\int_0^\infty \int_0^\infty
\Psi^{\alpha}_{q,g}(f)(a,b)g_{a,b}(x)b^{2\alpha+1}\frac{d_qad_qb}{a^2},~~~~~~~~~x\in\R_{q,+}.
\end {equation}
\end {theorem}
\proof For $x\in\R_{q,+}$, we have $h= \delta_x$ belongs to
$\L_q^2(\R_{q,+}, x^{2\alpha +1}d_qx)$. On the other hand,
according to the relation (\ref{Parseval}) of the previous
theorem,  the definition of $\Psi^{\alpha}_{q,g}$ and the
definition of the $q$-Jackson integral, we have
\begin {eqnarray*}
(1-q)x^{2\alpha+2}f(x)&=&\int_0^\infty
f(t)\overline{h}(t)t^{2\alpha+1}d_qt=\frac{1}{C_g}\int_0^\infty
\int_0^\infty \Psi^{\alpha}_{q,g}(f)
(a,b)\overline{\Psi^{\alpha}_{q,g}}(h)
(a,b)b^{2\alpha+1}\frac{d_qad_qb}{a^2}.\\&=&\frac{c_{\alpha,q}}{C_g}\int_0^\infty
\int_0^\infty \Psi^{\alpha}_{q,g}(f)(a,b)\left(\int_0^\infty
\overline{h}(t)g_{a,b}(t)t^{2\alpha+1}d_qt\right)b^{2\alpha+1}\frac{d_qad_qb}{a^2}\\
&=& (1-q)x^{2\alpha+2}\frac{c_{\alpha,q}}{C_g}\int_0^\infty
\int_0^\infty
\Psi^{\alpha}_{q,g}(f)(a,b)g_{a,b}(x)b^{2\alpha+1}\frac{d_qad_qb}{a^2}.
\end {eqnarray*}
Thus
$$ f(x)= \frac{c_{\alpha,q}}{C_g}\int_0^\infty \int_0^\infty
\Psi^{\alpha}_{q,g}(f)(a,b)g_{a,b}(x)b^{2\alpha+1}\frac{d_qad_qb}{a^2}.$$
Which completes the proof.\endproof
\section{Coherent states}
Theorem \ref{perplan} shows that the continuous wavelet transform
associated with the $q$-Bessel operator $\Psi^{\alpha}_{q,g}$ is
an isometry from the Hilbert space $\L_q^2(\R_{q,+}, x^{2\alpha
+1}d_qx)$ into the Hilbert space
$\L_q^2(\R_{q,+}\times\widetilde{\R}_{q,+};
b^{2\alpha+1}\frac{d_qad_qb}{a^2C_g})$ (the space of square
integrable functions on $\R_{q,+}\times\widetilde{\R}_{q,+}$ with
respect to the measure $b^{2\alpha+1}\frac{d_qad_qb}{a^2C_g}$).
For the characterization of the image of $\Psi^{\alpha}_{q,g}$, we
consider the vectors $g_{a,b}, ~~~~(a,b)\in
\R_{q,+}\times\widetilde{\R}_{q,+}$, as a set of coherent states
in the Hilbert space $\L_q^2(\R_{q,+}, x^{2\alpha +1}d_qx)$ (see
\cite{koornwinder}).
\begin{defin}\label{definition}
A set of coherent states in a Hilbert space $\mathcal{H}$ is a
subset $\{g_\textit{l}\}_{\textit{l}\in \mathcal{L}}$ of
$\mathcal{H}$ such that\\
\indent i) $\mathcal{L}$ is a locally  compact topological space
and the mapping $\textit{l}\mapsto g_\textit{l}$ is continuous
from $\mathcal{L}$ into $\mathcal{H}$.\\
\indent ii) There is a positive Borel measure $d\textit{l}$ on
$\mathcal{L}$ such that, for $f\in \mathcal{H}$,
$$ \parallel f\parallel^2=\int_\mathcal{L}\mid
(f,g_\textit{l})\mid^2d\textit{l},$$ where $(.,.)$ and $\parallel
. \parallel $ are respectively the scalar product and the norm of
$\mathcal{H}$.
\end{defin}

Let now $\mathcal{H}= \L_q^2(\R_{q,+}, x^{2\alpha +1}d_qx)$,
$\mathcal{L}=\R_{q,+}\times\widetilde{\R}_{q,+}$ equipped with the
induced topology of $\R^2$.\\
Choose a nonzero function $g\in \L_q^2(\R_{q,+}, x^{2\alpha
+1}d_qx)$ and let $g_\textit{l}=g_{a,b}$, $\textit{l}=(a,b)\in
\mathcal{L}$ be  given by the relation (\ref{wavelet}). Then we
have a set of coherent states. Indeed, i) of Definition
\ref{definition} is satisfied, because of Proposition \ref{cont},
and ii) of Definition \ref{definition} is satisfied, for the
measure $\ds b^{2\alpha+1}\frac{d_qad_qb}{a^2C_g}$ (see Theorem
\ref{perplan}). By adaptation of the approach introduced by T. H.
Koornwinder in \cite{koornwinder}, we obtain the following result:

\begin{theorem} Let $F$ be in
$\L_q^2(\R_{q,+}\times\widetilde{\R}_{q,+};
b^{2\alpha+1}\frac{d_qad_qb}{a^2C_g})$. Then $F$ belongs to $Im
\Psi^{\alpha}_{q,g}$ if and only if
\begin{equation}
F(a,b)=\frac{1}{C_g}\int_0^\infty \int_0^\infty
F(a',b')\left(\int_0^\infty
g_{a',b'}(x)\overline{g_{a,b}}(x)x^{2\alpha+1}d_qx\right)(b')^{2\alpha+1}\frac{d_qa'd_qb'}{(a')^2}.
\end{equation}
\end{theorem}

\section{Inversion formulas for the $q$-Riemann-Liouville and the $q$-Weyl operators}
{\bf Notations.} We denote by

\begin{itemize}
\item  $\mathcal{S}_{\ast q,\alpha }(\mathbb{R}_q)$ the subspace of $%
\mathcal{S}_{\ast q}(\mathbb{R}_q)$ constituted  of functions $f$
such that
\begin{equation*}
\int_{0}^{\infty }f\left( x\right) x^{2k+2\alpha +1}d_{q}x=0,\
k=0,1,...\ .
\end{equation*}

\item $\mathcal{S}_{\ast q}^{0}(\mathbb{R}_q)$ the subspace of $\mathcal{S}%
_{\ast q}(\mathbb{R}_q)$ constituted of functions $f$ such that%
\begin{equation*}
 D_{q}^{2k}f(0) =0,~~~~ k=0,1,...\ .
\end{equation*}
\end{itemize}

The q-Riemann-Liouville  transform $R_{\alpha,q}$ is defined on
$\mathcal{D}_{* q}({\R}_q)$ by (see \cite{Fih})
\begin {equation} R_{\alpha ,q}\left( f\right) (x)=
    \frac{(1+q)\Gamma
_{q^{2}}\left( \alpha +1\right) }{\Gamma _{q^{2}}\left(
\frac{1}{2}\right) \Gamma _{q^{2}}\left(\alpha
+\frac{1}{2}\right)}  \int_{0}^{1}\frac{\left(
t^{2}q^{2};q^{2}\right) _{\infty }}{\left( t^{2}q^{2\alpha
+1};q^{2}\right) _{\infty }}f(xt)d_{q}t.
\end{equation}
 The q-Weyl transform is defined on
$\mathcal{D}_{* q}({\R}_q)$ by (see \cite{Fih})
\begin {equation}
W_{\alpha ,q}\left( f\right) \left( x\right)
=\frac{q(1+q^{-1})^{-\alpha+\frac{1}{2}}\Gamma _{q^{2}}\left(
\alpha +1\right) }{ \Gamma _{q^{2}}^2\left(\alpha
+\frac{1}{2}\right)}\int_{qx}^{\infty }\frac{\left(
x^{2}/t^{2}q^{2};q^{2}\right) _{\infty }}{\left( q^{2\alpha
+1}x^{2}/t^{2};q^{2}\right) _{\infty }}f(t)t^{2\alpha }d_{q}t.
\end{equation}
These two operators are isomorphism on $\mathcal{D}_{* q}({\R}_q)$
and we have (see \cite{Fih}) $$ \Delta_{\alpha, q}\circ R_{\alpha
,q}=R_{\alpha ,q}\circ\Delta_{q}.$$
 and $$R_{\alpha ,q}(f*_qg)=R_{\alpha ,q}(f)*_BR_{\alpha
 ,q}(g),~~~~f,g \in \mathcal{D}_{*
q}({\R}_q),$$
 where "$*_q$" is the $q$-even convolution product associated with
 the operator  $\Delta_{q}$ studied in \cite{FB}.\\
 The $q$-Fourier-cosine transform $\F_q$
(studied in \cite{FB}) and the $q$-Bessel transform are linked by
 the following relation (see \cite{Fih}):
\begin{propo}\label{omar1}
For $f\in \mathcal{S}_{\ast q}({\R}_q),$ we have
\begin{eqnarray}
\f( f) &=&\F_q\circ W_{\alpha ,q}( f).
\end{eqnarray}
\end{propo}
We state the following results, useful in the sequel.
\begin{theorem}\label{w1}
The $q$-Fourier-cosine transform $\F_q$ is a topological
isomorphism from $\mathcal{S}_{\ast q,-1/2}(\mathbb{R}_q)$ into
$\mathcal{S}_{\ast q}^{0}(\mathbb{R}_q)$.
\end{theorem}
\proof From the Plancheral formula ( see \cite{D}), $\F_q$ is a
topological isomorphism from $\mathcal{S}_{\ast q}(\mathbb{R}_q)$
into itself.
 Moreover, using the fact that $D_q^2\cos(x;q^2)= -\cos(qx;q^2)$,
  one can prove by induction  that for $n\in\N$ and   $f$ in $\mathcal{S}_{\ast q}(\mathbb{R}_q)$,
 there exists a constant $C_{q,n}$, such that
 $$ D_q^{2n}\F_q(f)(0)= C_{q,n}\int_0^\infty f(t)t^{2n}d_qt,$$
 which achieves the proof.
\endproof
Similarly, we have the following result.
\begin{theorem}\label{w2}
The $q$-Fourier-Bessel transform $\f$ is a topological isomorphism
from $\mathcal{S}_{\ast q,\alpha}(\mathbb{R}_q)$ into
$\mathcal{S}_{\ast q}^{0}(\mathbb{R}_q)$.
\end{theorem}
\begin{coro}\label{w3}
The $q$-Weyl transform $W_{\alpha ,q}$ is a topological
isomorphism from $\mathcal{S}_{\ast q,\alpha}(\mathbb{R}_q)$ into
$\mathcal{S}_{\ast q,-1/2}(\mathbb{R}_q)$.
\end{coro}
\proof From the relation $\f=\F_q\circ W_{\alpha ,q}$, one can see
that $$W_{\alpha ,q}=\F_q^{-1}\circ\f.$$ We deduce the result from
this relation and Theorems \ref{w1} and \ref{w2}.
\endproof
\begin{propo}
For $f$ in $\mathcal{S}_{\ast q,-1/2}(\mathbb{R}_q)$ (resp.
$\mathcal{S}_{\ast q,\alpha}(\mathbb{R}_q)$) and $g$ in
$\mathcal{S}_{\ast q}(\mathbb{R}_q)$ the function $f*_qg$ (resp.
$f*_Bg$) belongs to $\mathcal{S}_{\ast q,-1/2}(\mathbb{R}_q)$
(resp. $\mathcal{S}_{\ast q,\alpha}(\mathbb{R}_q)$).
\end{propo}
\proof The proof follows from Theorem \ref{w1} (resp. \ref{w2})
and the fact that \\ $f*_qg =\F_q(\F_q(f).\F_q(g))$ (resp. $f*_Bg
=q^{-4\alpha-2}\f(\f(f).\f(g))$.
\endproof
\begin{propo}\label{ww1}
The operator $K_{\alpha,q,1 }$ defined by
\begin{equation*}
K_{\alpha ,q,1}( f)  =\frac{\Gamma _{q^2}\left( 1/2\right)
}{q^{3\alpha +3/2}(1+q)^{(\alpha +1/2)}\Gamma _{q^{2}}( \alpha +1)
}\F_q^{-1}( \left\vert \lambda \right\vert ^{2\alpha +1}\F_q( f) )
\end{equation*}
is a topological isomorphism from $\mathcal{S}_{\ast q,-1/2}(\mathbb{\R}_q%
) $ into itself.
\end{propo}
\proof The multiplication operator $$f\mapsto \frac{\Gamma
_{q^2}\left( 1/2\right) }{q^{3\alpha +3/2}(1+q)^{(\alpha
+1/2)}\Gamma _{q^{2}}( \alpha +1) } \left\vert \lambda \right\vert
^{2\alpha +1} f$$
  is a topological isomorphism from $\mathcal{S}_{\ast
  q}^0(\mathbb{\R}_q)$ into itself. The inverse is given by
  $$f\mapsto \frac{q^{3\alpha +3/2}(1+q)^{(\alpha +1/2)}\Gamma _{q^{2}}(
\alpha +1) }{\Gamma _{q^2}\left( 1/2\right)\left\vert \lambda
\right\vert ^{2\alpha +1} }  f.$$ The result follows from Theorem
\ref{w1}.\endproof
\begin{propo}\label{ww2}
The operator $K_{\alpha,q,2 }$ defined by
\begin{equation*}
K_{\alpha, q, 2}( f) ( x) =\frac{\Gamma _{q^{2}}\left( 1/2\right)
}{q^{3\alpha +3/2}(1+q)^{(\alpha +1/2)}\Gamma _{q^{2}}( \alpha +1)
}\f^{-1}( \left\vert \lambda \right\vert ^{2\alpha +1}\f( f) ) (x)
\end{equation*}%
is a topological isomorphism from $\mathcal{S}_{*q,\alpha
}(\mathbb{\R}_q)$ into itself.
\end{propo}
\proof From the relation $\f=\F_q\circ W_{\alpha,q}$ and  the
definition of $K_{\alpha,q,1 }$, we have for all $f\in
\mathcal{S}_{*q,\alpha }(\R_q)$
\begin {equation}
K_{\alpha,q,2 }=W_{\alpha,q}^{-1}\circ K_{\alpha,q,1 }\circ
W_{\alpha,q}.
\end{equation}
We deduce the result from Proposition \ref{ww1}  and Corollary
\ref{w3}.
\endproof
\begin{propo}\label{www1}
i)  For all $f\in \mathcal{S}_{* q,-1/2}(\R_q)$ and $g\in
\mathcal{S}_{* q}(\R_q),$ we have
\begin{equation*}
K_{\alpha ,q,1}(f* _{q}g) =K_{\alpha ,q,1}( f) * _{q}g.
\end{equation*}

ii)  For all $f\in \mathcal{S}_{* q,\alpha }(\R_q)$ and $g\in
\mathcal{S}_{* q}(\R_q)$, we have
\begin{equation*}
K_{\alpha ,q,2}( f*_Bg) =K_{\alpha ,q,2}( f) * _Bg.
\end{equation*}
\end{propo}

\proof It suffices to prove one of the two relations. We have
\begin{eqnarray*}
K_{\alpha ,q,1}\left( f* _{q}g\right) &=&\frac{\Gamma
_{q^{2}}\left( 1/2\right) }{q^{3\alpha +3/2}(1+q)^{(\alpha
+1/2)}\Gamma _{q^{2}}( \alpha +1) }\F_q^{-1}\left( \left\vert
\lambda \right\vert ^{2\alpha +1}\F_q\left( f* _{q}g\right)
\right)
\\
&=&\frac{\Gamma _{q^{2}}\left( 1/2\right) }{q^{3\alpha
+3/2}(1+q)^{(\alpha +1/2)}\Gamma _{q^{2}}( \alpha +1)
}\F_q^{-1}\left( \left\vert \lambda \right\vert ^{2\alpha
+1}\F_q\left( f\right)
\F_q\left( g\right) \right) \\
&=&\frac{\Gamma _{q^{2}}\left( 1/2\right) }{q^{3\alpha
+3/2}(1+q)^{(\alpha +1/2)}\Gamma _{q^{2}}( \alpha +1) }
 \left\{ \F_q^{-1}\left(\vert \lambda \vert^{2\alpha +1}\F_q( f)\right)
\right\} * _{q}g
 \\
&=&K_{\alpha ,q,1}( f)* _{q}g.
\end{eqnarray*}
\endproof

\begin{theorem}\label{www2}
For all $f\in \mathcal{S}_{\ast q,\alpha }({\R}_q),$ we have the
following inversion formulas for the operator $R_{\alpha,q }$
\begin {equation}\label{invr}
f=R_{\alpha ,q}\circ K_{\alpha ,q,1}\circ W_{\alpha ,q}( f)
\end{equation}
\begin {equation}\label{invr0}
f= R_{\alpha ,q}\circ W_{\alpha ,q}\circ K_{\alpha ,q,2}( f).
\end{equation}
\end{theorem}
\proof Using the properties of the operator $R_{\alpha ,q}$,
studied in \cite{Fih}, Theorem \ref{ld} and Proposition
\ref{omar1}, we obtain for $x\in\widetilde{\R}_{q,+}$,
\begin{eqnarray*}
q^{4\alpha+2}f(x)&=& c_{\alpha,q}\int_0^\infty \f (f)
(\lambda)\ji(\lambda
x;q^2)\lambda^{2\alpha+1}d_q\lambda\\
&=&R_{\alpha ,q}\left[c_{\alpha,q} \int_0^\infty \f (f)
(\lambda)\cos(\lambda
\centerdot;q^2)\lambda^{2\alpha+1}d_q\lambda\right](x)\\
&=&  R_{\alpha ,q}\left[c_{\alpha,q} \int_0^\infty
\lambda^{2\alpha+1}\F_q\circ W_{\alpha ,q}(f)
(\lambda)\cos(\lambda
\centerdot;q^2)d_q\lambda\right](x)\\
&=&  R_{\alpha ,q} \left\{\frac{c_{\alpha,q}}{c_{-1/2,
q}}\F_q^{-1} \left[\lambda^{2\alpha+1}\F_q\circ W_{\alpha ,q}(f)
\right]\right\}(x)\\
 &=&q^{4\alpha +2} R_{\alpha ,q}
\left\{\F_q^{-1}
\left[\frac{\Gamma_q(1/2)\lambda^{2\alpha+1}}{q^{3\alpha+3/2}(1+q)^{\alpha+1/2}\Gamma_q(\alpha+1)}\F_q\circ
W_{\alpha ,q}(f) \right]\right\}(x).
\end{eqnarray*}
Thus, $\forall x\in\widetilde{\R}_{q,+},~~~~f(x)= R_{\alpha
,q}\circ K_{\alpha ,q,1}\circ W_{\alpha ,q}( f)(x).$\\ We deduce
the second relation from the first relation and the fact
$$K_{\alpha ,q,2}=W_{\alpha ,q}^{-1}\circ K_{\alpha ,q,1}\circ
W_{\alpha ,q}.$$
\endproof
\begin{coro}\label{Zeineb}
The operator $R_{\alpha,q }$ is a topological isomorphism from
$\mathcal{S}_{* q,-1/2 }({\R}_q)$ into $\mathcal{S}_{* q,\alpha
}({\R}_q)$.
\end{coro}
\proof We deduce the result from Proposition \ref{ww1}, Corollary
\ref{w3} and the relation\\ (\ref{invr}).\endproof
 Similarly, we have the following result.
\begin{theorem} \label{nej}
For all $f\in \mathcal{S}_{\ast q,-1/2}({\R}_q),$  we have the
following inversion formulas for the operator $W_{\alpha ,q}$
\begin {equation}
f=W_{\alpha ,q}\circ R_{\alpha ,q}\circ K_{\alpha ,q,1}( f).
\end{equation}
\begin {equation}
f=W_{\alpha ,q}\circ K_{\alpha ,q,2}\circ R_{\alpha ,q}( f)
\end {equation}

\end{theorem}
\proof For $f\in \mathcal{S}_{\ast q,-1/2}({\R}_q)$, Corollary
\ref{w3} (resp. (\ref{Zeineb})) implies that $W_{\alpha ,q}^{-1}(
f)$ (resp. $R_{\alpha ,q}( f)$) belongs to $\mathcal{S}_{\ast
q,\alpha}({\R}_q).$ Then by writing the relation (\ref{invr}) (
resp. \ref{invr0})  for $W_{\alpha ,q}^{-1}( f)$ (resp. $R_{\alpha
,q}( f)$), we obtain the result.
\endproof

\begin{coro}\label{299}
i)For all $f,g\in \mathcal{S}_{* q,\alpha}(\R_q)$, we have
\begin {equation}\label{omar2}
W_{\alpha ,q}\left( f* _Bg\right) =W_{\alpha ,q}( f)
* _{q}W_{\alpha ,q}( g).
\end{equation}

ii) For all $f, g\in \mathcal{S}_{* q, -1/2}(\R_q)$  we have
\begin {equation}\label{omar3}
R_{\alpha ,q}\left( f* _{q}g\right) = R_{\alpha ,q}\left( f\right)
* _BW_{\alpha ,q}^{ -1}( g).
\end{equation}
\end{coro}

\proof i) From Proposition \ref{omar1}, we have
\begin{eqnarray*}
W_{\alpha ,q}\left( f* _Bg\right) &=&\F_q^{-1}\circ
\f( f* _Bg) \\
&=&\F_q^{-1}\left( \f( f) \f( g)\right)\\ &=&\F_q^{-1}\circ
\f( f) * _{q}\F_q^{-1}\f( g) \\
&=&W_{\alpha ,q}( f) * _{q}W_{\alpha ,q}( g).
\end{eqnarray*}%
ii) Using Theorem \ref{nej} and Proposition \ref{www1}, we obtain
\begin{eqnarray*}
R_{\alpha ,q}^{-1}\left( R_{\alpha ,q}( f) *_B W_{\alpha ,q}^{
-1}( g) \right) &=&W_{\alpha ,q}\circ K_{\alpha ,q,2}\left(
R_{\alpha ,q}( f) *_B
W_{\alpha ,q}^{-1}( g) \right) \\
&=&W_{\alpha ,q}\left( K_{\alpha ,q,2}\circ R_{\alpha ,q}(
f) *_BW_{\alpha ,q}^{ -1}( g) \right) \\
&=&W_{\alpha ,q}\circ K_{\alpha ,q,2}\circ R_{\alpha ,q}( f) *
_{q} g.
\end{eqnarray*}%
On the other hand, we have
\begin{equation*}
W_{\alpha ,q}\circ K_{\alpha ,q,2}\circ R_{\alpha ,q}( f) =f.
\end{equation*}%
So,
\begin{equation*}
R_{\alpha ,q}^{-1}\left( R_{\alpha ,q}( f) *_B W_{\alpha ,q}^{\
-1}( g) \right) =f* _{q}g.
\end{equation*}%
This achieves  the proof.
\endproof
\section{Inversion formulas for the $q$-Riemann-Liouville and the $q$-Weyl operators using wavelets }
We recall that (see \cite{FN}):\\
$\bullet$ the dilatation operator is defined for $a\in\R_{q,+}$ by
\begin {equation}
H_a(f)(x)=\frac{1}{\sqrt{a}}f\left(\frac{x}{a}\right),
\end{equation}
$\bullet$ a $q$-wavelet is an even and square $q$-integrable
function $g$ satisfying $$ 0<C_g^c=\int_0^\infty
|\F_q(g)|^2(a)\frac{d_qa}{a}<\infty,$$
 $\bullet$ for a $q$-wavelet $g$, the continuous $q$-wavelet transform ({\bf ie.} associated with the
operator $\Delta_q$) is defined on
$\R_{q,+}\times\widetilde{\R}_{q,+}$ by
\begin{eqnarray}
\Phi_{q,g}(f)(a,b)&=&c_{-1/2, q} \int_0^\infty
f(x)\overline{g_{a,b}^c}(x)d_qx\\\label{*}
&=&f*_q\overline{H_a(g)}(b),
\end{eqnarray}
where $g_{a,b}^c=T_{q,b}(H_a(g))$ and $T_{q,b}=T_{q,b}^{-1/2}$ is
the $q$-even translation operator studied in \cite{FB}.
\begin{propo}
For all $a\in\R_{q,+}$ and $g\in
L_q^2(\R_{q,+},x^{2\alpha+1}d_qx)$, we have\\
 1) $\ds g_a= \frac{1}{a^{2\alpha +3/2}}H_a(g)$;\\
 2) \begin {eqnarray}
  g_a &=&\frac{q^{-4\alpha-2}}{\sqrt{a}}\f \circ H_{a^{-1}}\circ
  \f(g)\\\label{*1}
  &=&  \frac{1}{\sqrt{a}}W_{\alpha, q}^{-1}\circ H_a\circ W_{\alpha,
  q}(g).
    \end{eqnarray}
\end{propo}
\proof 1) is clear.\\
2) From the facts $ \f(g_a)(\lambda)=\f(g)(a\lambda)$ and
$\F_q\circ H_a =H_{a^{-1}} \circ\F_q$ (see \cite{FN}), one can see
$$\f(g_a)= \frac{1}{\sqrt{a}}H_{a^{-1}} \circ\f(g).$$
Then, by using Proposition \ref{omar1} and Theorem \ref{ld}, we
obtain
\begin {eqnarray*}
g_a&=&\frac{1}{\sqrt{a}}\f^{-1}\circ H_{a^{-1}} \circ\f(g)=
\frac{q^{-4\alpha-2}}{\sqrt{a}}\f \circ H_{a^{-1}}\circ
  \f(g)\\
&=&\frac{1}{\sqrt{a}}W_{\alpha, q}^{-1}\circ H_a\circ W_{\alpha,
  q}(g).
\end  {eqnarray*}
\endproof
\begin{propo}
Let $g$ be a $q$-wavelet associated with the $q$-Bessel operator
in $\mathcal{S}_{*q, \alpha}(\R_q)$. Then for all $f$ in
$\mathcal{S}_{*q, \alpha}(\R_q)$, we have the following relation
\begin {equation}
\Psi_{q,g}^\alpha
(f)(a,.)=W_{\alpha,q}^{-1}\left[\Phi_{q,W_{\alpha,q}(g)}\left(W_{\alpha,q}(f)\right)(a,.)\right],
\qquad a\in \R_{q,+} .
\end{equation}
\end{propo}
\proof Let $a\in \R_{q,+}$, we have from the relations
(\ref{omar2}), (\ref{*}) and (\ref{*1}),
\begin{eqnarray*}
\Psi_{q,g}^\alpha
(f)(a,.)&=&\sqrt{a}f*_B\overline{g_a}=\sqrt{a}W_{\alpha,q}^{-1}\left[W_{\alpha,q}(f)*_qW_{\alpha,q}(\overline{g_a})
\right]\\&=&
W_{\alpha,q}^{-1}\left[W_{\alpha,q}(f)*_q\overline{H_a\circ
W_{\alpha,q}(g)} \right]\\
&=&W_{\alpha,q}^{-1}\left[\Phi_{q,W_{\alpha,q}(g)}\left(W_{\alpha,q}(f)\right)(a,.)\right].
\end{eqnarray*}
\endproof
\begin{theorem}
Let $g$ be a $q$-wavelet associated with the $q$-Bessel operator
in $\mathcal{S}_{*q, \alpha}(\R_q)$. Then \\
 1) for all $f$ in $\mathcal{S}_{*q, \alpha}(\R_q)$, we have the following relation
\begin {equation}
\Psi_{q,g}^\alpha
(f)(a,b)=R_{\alpha,q}\left[\Phi_{q,W_{\alpha,q}(g)}\left(R_{\alpha,q}^{-1}(f)\right)(a,.)\right](b),
\qquad a\in \R_{q,+},~~b\in \widetilde{\R}_{q,+};
\end{equation}

2) for all $f$ in $\mathcal{S}_{*q, -1/2}(\R_q)$, we have
\begin {equation}
\Phi_{q,W_{\alpha,q}(g)}
(f)(a,b)=W_{\alpha,q}\left[\Psi_{q,g}^\alpha\left(W_{\alpha,q}^{-1}(f)\right)(a,.)\right](b),
\qquad a\in \R_{q,+},~~b\in \widetilde{\R}_{q,+}.
\end{equation}
\end{theorem}
\proof 1)  From Corollary \ref{299} and the relations (\ref{*})
and (\ref{*1}), we obtain for $a\in \R_{q,+}$ and  $b\in
\widetilde{\R}_{q,+}$,
\begin{eqnarray*}
\Psi_{q,g}^\alpha (f)(a,b)&=&\sqrt{a}f*_B\overline{g_a}(b)\\
&=& \sqrt{a} R_{\alpha,q}\left[R_{\alpha,q}^{-1}(f)*_q
W_{\alpha,q}(\overline{g_a})\right](b)\\
&=&R_{\alpha,q}\left[R_{\alpha,q}^{-1}(f)*_q
\overline{H_a\circ W_{\alpha,q}(g)}\right](b)\\
&=&
R_{\alpha,q}\left[\Phi_{q,W_{\alpha,q}(g)}\left(R_{\alpha,q}^{-1}(f)\right)(a,.)\right](b).
\end{eqnarray*}
2) For $a\in \R_{q,+}$ and  $b\in \widetilde{\R}_{q,+}$, we have
by the relation (\ref{*}), Corollary \ref{299} and the relation
(\ref{*1}),
\begin{eqnarray*}
\Phi_{q,W_{\alpha,q}(g)}(f)(a,b)&=& f*_q \overline{H_a\circ
W_{\alpha,q}(g)}(b)\\
&=& W_{\alpha,q}\left[W_{\alpha,q}^{-1}(f)*_B
\overline{W_{\alpha,q}^{-1}\circ H_a\circ
W_{\alpha,q}(g)}\right](b)\\
&=&W_{\alpha,q}\left[\sqrt{a}W_{\alpha,q}^{-1}(f)*_B
\overline{g_a}\right](b)\\
&=&W_{\alpha,q}\left[\Psi_{q,g}^\alpha\left(W_{\alpha,q}^{-1}(f)\right)(a,.)\right](b).
\end{eqnarray*}
\endproof
\begin{propo}\label{ne1} 1) If $g$ is a  $q$-wavelet in $\mathcal{S}_{*q,
-1/2}(\R_q)$, then $K_{\alpha, q, 1}(g)$ is  a  $q$-wavelet in
$\mathcal{S}_{*q, -1/2}(\R_q)$ and we have
\begin {equation}
K_{\alpha, q, 1}\circ H_a(g)= \frac{1}{a^{2\alpha+1}}H_a\circ
K_{\alpha, q, 1}(g),~~~~a\in\R_{q,+}.
\end{equation}
2) If $g$ is a  $q$-wavelet associated with the $q$-Bessel
operator in $\mathcal{S}_{*q, \alpha}(\R_q)$, then $K_{\alpha, q,
2}(g)$ is  a $q$-wavelet in $\mathcal{S}_{*q, \alpha}(\R_q)$ and
we have
\begin {equation}
K_{\alpha, q, 2}(g_a)= \frac{1}{a^{2\alpha+1}}( K_{\alpha, q,
2}(g))_a,~~~~a\in\R_{q,+}.
\end{equation}
\end{propo}
\proof 1) Let $g$ be a  $q$-wavelet in $\mathcal{S}_{*q,
-1/2}(\R_q)$. From the definition of $K_{\alpha, q, 1}$, we have
for $\lambda\in\R_{q,+}$,
$$ \F_q(K_{\alpha, q, 1}(g))(\lambda)=\frac{\Gamma _{q^2}\left( 1/2\right)
}{q^{3\alpha +3/2}(1+q)^{(\alpha +1/2)}\Gamma _{q^{2}}( \alpha +1)
} \lambda^{2\alpha +1}\F_q( g)(\lambda).$$ Proposition 4 of
\cite{FN}, implies that $K_{\alpha, q, 1}(g)$ is a $q$-wavelet. On
the other hand, using the fact $\F_q\circ H_a=
H_{a^{-1}}\circ\F_q,~~ a\in\R_{q,+}$ and the above equality, we
obtain
\begin {eqnarray*}
\F_q(H_a\circ K_{\alpha, q, 1}(g))(\lambda)&=& \frac{\Gamma
_{q^2}\left( 1/2\right) }{q^{3\alpha +3/2}(1+q)^{(\alpha
+1/2)}\Gamma _{q^{2}}( \alpha +1)
}H_{a^{-1}}\left(\lambda^{2\alpha +1}\F_q( g)(\lambda)\right)\\
&=& a^{2\alpha+1}\frac{\Gamma _{q^2}\left( 1/2\right) }{q^{3\alpha
+3/2}(1+q)^{(\alpha +1/2)}\Gamma _{q^{2}}( \alpha +1)
}\lambda^{2\alpha +1}H_{a^{-1}}\left(\F_q( g)(\lambda)\right)\\
&=& a^{2\alpha+1}\frac{\Gamma _{q^2}\left( 1/2\right) }{q^{3\alpha
+3/2}(1+q)^{(\alpha +1/2)}\Gamma _{q^{2}}( \alpha +1)
}\lambda^{2\alpha +1}\F_q\left(H_a (g)\right)(\lambda),
\end{eqnarray*}
which gives the result.\\
 2)  Let $g$ be a  $q$-wavelet associated with the $q$-Bessel operator in $\mathcal{S}_{*q,
\alpha}(\R_q)$. From the definition of $K_{\alpha, q, 2}$, we have
for $\lambda\in\R_{q,+}$,
$$ \f(K_{\alpha, q, 2}(g))(\lambda)=\frac{\Gamma _{q^2}\left( 1/2\right)
}{q^{3\alpha +3/2}(1+q)^{(\alpha +1/2)}\Gamma _{q^{2}}( \alpha +1)
} \lambda^{2\alpha +1}\f( g)(\lambda).$$ Proposition \ref{ne},
implies that $K_{\alpha, q, 2}(g)$ is a $q$-wavelet associated
with the $q$-Bessel operator.\\
Moreover, for $\lambda\in\R_{q,+}$, we have
\begin {eqnarray*}
\f(K_{\alpha, q, 2}(g_a))(\lambda)&=&\frac{\Gamma _{q^2}\left(
1/2\right) }{q^{3\alpha +3/2}(1+q)^{(\alpha +1/2)}\Gamma _{q^{2}}(
\alpha +1) } \lambda^{2\alpha +1}\f( g_a)(\lambda)\\
&=&\frac{\Gamma _{q^2}\left( 1/2\right) }{q^{3\alpha
+3/2}(1+q)^{(\alpha +1/2)}\Gamma _{q^{2}}(
\alpha +1) } \lambda^{2\alpha +1}\f( g)(a\lambda)\\
&=&\frac{1}{a^{2\alpha+1}}\f([K_{\alpha, q, 2}(g)]_a)(\lambda).
\end{eqnarray*}
This achieves the proof.
\endproof
\begin{theorem}
Let $g$ be a $q$-wavelet  associated with the $q$-Bessel
operator in  $\mathcal{S}_{*q, \alpha}(\R_q)$. Then for $a\in\R_{q,+}$ and $b\in\widetilde{\R}_{q,+}$, we have:\\
 1) for all $f$ in $\mathcal{S}_{*q, \alpha}(\R_q)$,
 \begin {equation}
 \Psi_{q,g}^\alpha(f)(a,b)=\frac{1}{a^{2\alpha+1}}R_{\alpha,q}\left[
  \Phi_{q,K_{\alpha, q, 1}\circ W_{\alpha,q}(g)}(W_{\alpha,q}(f))(a,.)\right](b);
 \end{equation}

 2) for all $f$ in $\mathcal{S}_{*q, -1/2}(\R_q)$,
 \begin {equation}
 \Phi_{q,W_{\alpha,q}(g)}(f)(a,b)=\frac{1}{a^{2\alpha+1}}W_{\alpha,q}\left[
  \Psi_{q,K_{\alpha, q, 2}(g)}^\alpha(R_{\alpha,q}(f))(a,.)\right](b).
 \end{equation}
\end{theorem}
\proof 1) Let $f$ be  in $\mathcal{S}_{*q, \alpha}(\R_q)$,
$a\in\R_{q,+}$ and $b\in\widetilde{\R}_{q,+}$. Using Corollary
\ref{299}, we obtain
\begin {eqnarray*}
\Psi_{q,g}^\alpha(f)(a,b)&=&\sqrt{a}f*_B\overline{g_a}(b)\\
&=&
\sqrt{a}R_{\alpha,q}\left[W_{\alpha,q}(f)*_qR_{\alpha,q}^{-1}(\overline{g_a})\right](b).
\end{eqnarray*}
But, from Theorem \ref{www2}, Proposition \ref{ne1} and the
relation (\ref{*1}), we have
\begin {eqnarray*}
R_{\alpha,q}^{-1}(\overline{g_a})&=&
K_{\alpha, q, 1}\circ W_{\alpha,q}(\overline{g_a})\\
&=&\frac{1}{\sqrt{a}}K_{\alpha, q, 1}\circ \overline{H_a\circ W_{\alpha,q}(g)}\\
&=&\frac{1}{a^{2\alpha+3/2}}\overline{H_a\circ K_{\alpha, q,
1}\circ W_{\alpha,q}(g)}.
\end{eqnarray*}
Thus,
\begin {eqnarray*}
\Psi_{q,g}^\alpha(f)(a,b)&=&\frac{1}{a^{2\alpha+1}}R_{\alpha,q}\left[W_{\alpha,q}(f)*_q\overline{H_a\circ
K_{\alpha, q, 1}\circ W_{\alpha,q}(g)}\right](b)\\
&=&\frac{1}{a^{2\alpha+1}}R_{\alpha,q}\left[
  \Phi_{q,K_{\alpha, q, 1}\circ W_{\alpha,q}(g)}(W_{\alpha,q}(f))(a,.)\right](b).
\end{eqnarray*}
2) Let $f$ be  in $\mathcal{S}_{*q, -1/2}(\R_q)$, $a\in\R_{q,+}$
and $b\in\widetilde{\R}_{q,+}$. Using Corollary \ref{299}, we
obtain
\begin {eqnarray*}
\Phi_{q,W_{\alpha,q}(g)}(f)(a,b)&=&f*_q
\overline{H_a\circ W_{\alpha,q}(g)}\\
&=&W_{\alpha,q}[W_{\alpha,q}^{-1}(f)*_B
\overline{W_{\alpha,q}^{-1}\circ H_a\circ W_{\alpha,q}(g)}]\\
&=& \sqrt{a}  W_{\alpha,q}[W_{\alpha,q}^{-1}(f)*_B
\overline{g_a}]\\
&=& \sqrt{a}  W_{\alpha,q}[K_{\alpha, q, 2}\circ
R_{\alpha,q}(f)*_B \overline{g_a}],
\end{eqnarray*}
since $W_{\alpha,q}^{-1}(f)=K_{\alpha, q, 2}\circ R_{\alpha,q}(f)$
according to Theorem \ref{nej}.\\ Using Propositions \ref{www1}
and \ref{ne1}, we obtain
\begin {eqnarray*}
\Phi_{q,W_{\alpha,q}(g)}(f)(a,b)&=& \sqrt{a}
W_{\alpha,q}[R_{\alpha,q}(f)*_B \overline{K_{\alpha, q, 2}(g_a)}](b)\\
&=& \frac{1}{a^{2\alpha+1/2}} W_{\alpha,q}[R_{\alpha,q}(f)*_B
\overline{[K_{\alpha, q, 2}(g)]_a}](b)\\
&=& \frac{1}{a^{2\alpha+1}} W_{\alpha,q}[\Psi_{q,K_{\alpha, q,
2}(g)}^\alpha(R_{\alpha,q}(f))(a,.)](b)
\end{eqnarray*}
\endproof
\begin{theorem}
Let $g$ be a $q$-wavelet  associated with the $q$-Bessel operator
in  $\mathcal{S}_{*q, \alpha}(\R_q)$ and $x\in\R_{q,+}$. Then,\\
1) for all $f$ in $\mathcal{S}_{*q, -1/2}(\R_q)$, we have
\begin {eqnarray*}
W_{\alpha,q}^{-1}(f)(x)=\frac{c_{\alpha,q}}{C_g}\int_0^\infty\left(\int_0^\infty
R_{\alpha,q}[\Phi_{q,K_{\alpha, q, 1}\circ
W_{\alpha,q}(g)}(f)(a,.)](b)\times g_{a,b}(x)\frac{b^{2\alpha
+1}}{a^{2\alpha+3}}  d_qb\right)d_qa;
\end {eqnarray*}
2) for all $f$ in $\mathcal{S}_{*q, \alpha}(\R_q)$, we have
\begin {eqnarray*}
R_{\alpha,q}^{-1}(f)(x) = \frac{c_{-\frac{1}{2},
q}}{C^c_g}\int_0^\infty\left(\int_0^\infty
W_{\alpha,q}\left[\Psi_{q,K_{\alpha, q, 2}(g)}^\alpha (f)(a,.)
\right](b)g_{a,b}^c(x) \frac{d_qb}{a^{2\alpha+3}}\right) d_qa.
\end {eqnarray*}
\end{theorem}
\proof 1) is a simple deduction from the previous theorem and
Theorem \ref{omri}.\\
2) Similarly, the result derives from the previous theorem and
Theorem 7 of \cite{FN}.
\endproof

\end{document}